\input amstex
\documentstyle{amsppt}
\magnification\magstep1
\NoBlackBoxes 

\def\obr{^{-1}}
\def\sbs{\subset}

\def\Aut{\operatorname{Aut}}

\def\Dom{\operatorname{Dom}}
\def\Ran{\operatorname{Ran}}

\def\o{\omega}
\def\G{\Gamma}

\def\a{\alpha}
\def\d{\delta}
\def\D{\Delta}

\def\g{\gamma}

\def\l{\lambda}

\def\Th{\Theta}

\define\cl#1{\overline {#1}}

\def\sF{{\Cal F}}
\def\sG{{\Cal G}}

\def\sL{{\Cal L}}

\def\sN{{\Cal N}}

\def\sR{{\Cal R}}

\def\sU{{\Cal U}}
\def\sV{{\Cal V}}

\def\ti{\times}

\hyphenation{homo-ge-ne-ous}

\newcount\secno
\secno=-1
\newcount\subsecno
\subsecno=-1
\newcount\thmno
\thmno=0
\newcount\bigthmno
\bigthmno=0
\newcount\secthmno
\bigthmno=0

\def\ifundefined#1{\expandafter\ifx\csname#1\endcsname\relax}

\def\privateName#1{\ifundefined{privateName{#1}}\errmessage{Undefined name
#1}\else\csname privateName{#1}\endcsname\fi}

\def\name#1{\privateName{#1}}

\def\resetName#1#2{{\global\expandafter\edef\csname
privateName{#1}\endcsname{#2}}}

\def\setName#1#2{\ifundefined{privateName{#1}}\resetName{#1}{#2}\else
\errmessage{Name #1 is already in use}\fi}

\def\theoremName#1{%
\global\advance\thmno by1%
\global\advance\bigthmno by1%
\global\advance\secthmno by1%
\setName{#1}{%
        \ifnum\number\the\secno<0\else
               \ifnum\number\the\subsecno<0\else\number\the\subsecno.\fi
               \number\the\secno.%
        \fi
        \number\the\thmno
        }%
}

\def\theoremLabel#1{setName{#1}{#1}}

\def\No#1{\theoremName{#1}\name{#1}}

\def\label#1{\theoremLabel{#1}\privateName{#1}}

\def\newSection{\ifnum\secno<0%
                        \secno=1%
                \else
                        \advance\secno by1%
                \fi
                \ifnum\subsecno<0\else
                \subsecno=1%
                \fi
                \thmno=0%
                \secthmno=0}

\def\newSubSection{\ifnum\subsecno<0\subsecno=1\else
\advance\subsecno by1\fi
\thmno=0}

\def\setTheoremNumber#1{\thmno=#1\advance\thmno by-1}

\def\showTotalTheoremNumber{\message{{\number\the\bigthmno} claims etc.
total}}

\def\showSectionTheoremNumber{\message{The section contains
{\number\the\bigthmno} claims etc.}}


\def\setBibl#1#2{\setName{*bibl*#1}{#2}}
\def\bibl#1{\name{*bibl*#1}}

\topmatter
\author V. V. Uspenskij
\endauthor


\address 
Department of Mathematics, 321 Morton Hall, Ohio University, Athens OH~45701, USA
\endaddress
\email 
uspensk\@bing.math.ohiou.edu
\endemail



\title
The Roelcke compactification of groups of homeomorphisms
\endtitle

\abstract
Let $X$ be a zero-dimensional compact space such that all non-empty
clopen subsets of $X$ are homeomorphic to each other, and let $\Aut X$
be the group of all self-homeomorphisms of $X$ with the 
compact-open topology. We prove that the Roelcke compactification
of $\Aut X$ can be identified with the semigroup of all closed relations
on $X$ whose domain and range are equal to $X$. 
We use this to prove 
that the group $\Aut X$ is topologically simple and minimal, in the sense that it does not admit a strictly coarser
Hausdorff group topology. For $X=2^\o$ the last result is due to D.~Gamarnik.
\endabstract

\endtopmatter

\setBibl{dik}{1}
\setBibl{DPS}{2}
\setBibl{gamarn}{3}
\setBibl{RD}{4}
\setBibl{St}{5}
\setBibl{stoyan}{6}
\setBibl{minim}{7}

\document

\newSection

\head
\S~\number\the\secno. Introduction
\endhead

Let $G$ be a topological group. There are at least four natural uniform
structures on $G$ which are compatible with the topology [\bibl{RD}]: the left
uniformity $\sL$, the right uniformity $\sR$, their least upper bound
$\sL\vee\sR$ and their greatest lower bound $\sL\wedge\sR$. In [\bibl{RD}]
the uniformity $\sL\wedge\sR$ is called the {\it lower uniformity\/} on $G$;
we shall call it the {\it Roelcke uniformity\/}, as in [\bibl{stoyan}].
Let $\sN(G)$ be the filter of neighbourhoods of unity in $G$.
When $U$ runs over $\sN(G)$, the covers of the form
$\{xU:x\in G\}$, $\{Ux:x\in G\}$, $\{xU\cap Ux:x\in G\}$ and
$\{UxU:x\in G\}$ are uniform for $\sL$, $\sR$, $\sL\vee\sR$ and
$\sL\wedge\sR$, respectively, and generate the corresponding uniformity.

All topological groups  are assumed  to be Hausdorff. 
A  uniform  space  $X$  is  {\it
precompact\/} if its completion is  compact  or,  equivalently,  if every
uniform cover of $X$ has a finite subcover.
For any topological group $G$ the following are equivalent: (1) $G$ is 
$\sL$-precompact; (2) $G$ is $\sR$-precompact; (3) $G$ is 
$\sL\vee\sR$-precompact; (4) $G$ is a topological subgroup of a compact group.
If these conditions are satisfied, $G$ is said to be {\it precompact}.
Let us say that $G$ is {\it Roelcke-precompact\/} if $G$ is precompact with
respect to the Roelcke uniformity. A group $G$ is precompact if and only if
for every $U\in \sN(G)$ there exists a finite set $F\sbs G$ such that
$UF=FU=G$.
A group $G$ is Roelcke-precompact if and only if for every $U\in \sN(G)$ there
exists a finite $F\sbs G$ such that $UFU=G$.
Every precompact group is
Roelcke-precompact, but not vice versa.
For  example,  the  unitary  group  of  a  Hilbert  space  or  the group
$\operatorname{Symm}(E)$ of all permutations of a discrete set $E$,  both  with
the  pointwise convergence topology, are Roelcke-precompact but not precompact
[\bibl{stoyan}], [\bibl{RD}].
Unlike the usual precompactness, the property of being Roelcke-precompact is
not inherited by subgroups. 
(If $H$ is a subgroup of $G$, in general the Roelcke uniformity of $H$ is finer
than the uniformity induced on $H$ by the Roelcke uniformity of $G$.)
Moreover, every topological group is a subgroup
of a Roelcke-precompact group [\bibl{minim}].

The {\it Roelcke completion\/} of a topological group $G$ is the completion
of the uniform space $(G,\sL\wedge\sR)$. If $G$ is Roelcke-precompact,
the Roelcke completion of $G$ will be called the
{\it Roelcke compactification\/} of $G$.

A topological group is {\it minimal\/} if it does not
admit a strictly coarser Hausdorff group topology. 
Let  us  say  that  a  group $G$ is {\it topologically simple\/} if $G$ has no
closed normal subgroups besides $G$ and $\{1\}$. It was shown in 
[\bibl{stoyan}], [\bibl{minim}] that the Roelcke compactification of
some important topological groups has a natural structure of an ordered
semigroup with an involution, and that the study of this structure can be
used to prove that a given group is minimal and topologically simple.
In the present paper we apply this method to some groups of homeomorphisms.

A {\it semigroup\/} is a set with an associative binary operation.
Let $S$ be a semigroup with the multiplication $(x,y)\mapsto xy$.
We say that a self-map
$x\mapsto x^*$ of $S$ is an {\it involution\/} if $x^{**}=x$ and
$(xy)^*=y^*x^*$ for all $x,y\in S$. Every group has a natural involution
$x\mapsto x\obr$.
An element $x\in S$ is {\it symmetrical\/} if $x^*=x$, and a subset $A\sbs S$
is {\it symmetrical\/} if $A^*=A$.
An {\it ordered semigroup\/} is a semigroup with a partial order $\le$ such
that the conditions $x\le x'$ and $y\le y'$ imply $xy\le x'y'$.
An element $x\in S$ is {\it idempotent\/} if $x^2=x$.

Let $K$ be a compact space. A {\it closed relation\/} on $K$
is a closed subset of $K^2$.
Let $E(K)$ be the compact space  of all closed relations on $K$,
equipped with the Vietoris topology. 
If $R,S\in E(K)$, then the {\it composition\/} of $R$ and $S$ is
the relation $RS=\{(x,y):\exists z((x,z)\in S\text{ and }(z,y)\in R)\}$.
The relation $RS$ is closed, since it is the image of the closed
subset $\{(x,z,y): (x,z)\in S,\,(z,y)\in R\}$ of $K^3$ under the projection
$K^3\to K^2$ which is a closed map. 
If $R\in E(K)$, then the {\it inverse relation\/} $\{(x,y):(y,x)\in R\}$
will be denoted by $R^*$ or by $R\obr$; we prefer the first notation, since
we are interested in the algebraic structure on $E(K)$, and
in general $R\obr$ is not an inverse of $R$ in the algebraic sense.
The set $E(K)$ has a natural partial order. 
Thus $E(K)$ is an ordered semigroup
with an involution. In general the map $(R,S)\mapsto RS$ from
$E(K)^2$ to $E(K)$ is not (even separately) continuous.

For $R\in E(K)$ let $\Dom R=\{x: \exists y((x,y)\in R)\}$ and 
$\Ran R=\{y: \exists x((x,y)\in R)\}$. 
Put $E_0(K)=\{R\in E(K):\Dom R=\Ran R=K\}$. The set $E_0(K)$ is a closed
symmetrical subsemigroup of $E(K)$.

Denote by $\Aut(K)$
the group of all
self-homeomorphisms of $K$, equipped with the compact-open topology.
For  every  $h\in  \Aut(K)$ let $\G(h)=\{(x, h(x)): x\in K\}$ be the graph of
$h$.
The map $h\mapsto \G(h)$ from $\Aut(K)$ to $E_0(K)$ is a homeomorphic
embedding and a morphism of semigroups with an involution. 
Identifying every self-homeomorphism of $K$ with its graph, we consider
the group $\Aut(K)$ as a subspace of $E_0(K)$.  

We say that a compact space $X$ is {\it $h$-homogeneous}
if $X$ is zero-dimensional 
and all non-empty clopen subsets of $X$ are homeomorphic to each other.

\proclaim{\No{main}. Main Theorem}
Let $X$ be an $h$-homogeneous compact space,
and let $G=\Aut(X)$ be the topological group of all self-homeomorphisms
of $X$.
Then:
\roster
\item $G$ is Roelcke-precompact; 
the Roelcke compactification of $G$ can be identified with
the semigroup $E_0(X)$ of all closed relations $R$ on $X$ such that 
$\Dom R=\Ran R=X$;
\item  $G$ is minimal and topologically simple.
\endroster
\endproclaim

In the case when $X=2^\o$ is the Cantor set, 
the minimality of $\Aut(X)$ was proved by D.~Gamarnik [\bibl{gamarn}].

Let us explain how to deduce the second part of Theorem~\name{main} from
the first. 
Let $G=\Aut(X)$ be such as in 
Theorem~\name{main}, and let $f:G\to G'$ be a continuous onto homomorphism.
We must prove that either
$f$ is a topological isomorphism or $|G'|=1$. 
Let $\Th=E_0(X)$.
The first part of Theorem~\name{main} implies that 
$f$ can be extended
to a map $F:\Th\to \Th'$, where $\Th'$ is the Roelcke compactification of 
$G'$. 
Let $e'$ be the unity of $G'$, and     
let $S=F\obr(e')$. Then $S$ is a closed symmetrical subsemigroup of
$\Th$. 
Let $\D$ be the diagonal in $X^2$.
The set $\{r\in S:\D\sbs r\}$ has a largest element. Denote this element
by $p$. Then $p$ is a symmetrical idempotent in $\Th$ and hence 
an equivalence relation on $X$. The semigroup $S$ is invariant under
inner automorphisms of $\Th$, and so is the relation
$p$. But there are only two $G$-invariant closed
equivalence relations on $X$, namely
$\D$ and $X^2$. If $p=\D$, then $S\sbs G$, $G=F\obr(G')$ and $f$ is perfect.
Since $G$ has no non-trivial compact normal subgroups, we conclude that $f$
is a homeomorphism. If $p=X^2$, then $S=\Th$ and $G'=\{e'\}$.

A similar argument was used in [\bibl{minim}] 
to prove that every topological group
is a subgroup of a Roelcke-precompact topologically simple minimal group,
and in [\bibl{stoyan}] to yield an alternative proof of Stoyanov's theorem
asserting that the unitary group of a Hilbert space is minimal
[\bibl{St}], [\bibl{DPS}]. For more
information on minimal groups, see the recent survey by D.~Dikranjan 
[\bibl{dik}].

We prove the first part of Theorem~\name{main} in Section~2, and the second
part in Section~4.

\newSection

\head
\S~\number\the\secno. Proof of Main Theorem, part 1 
\endhead

Let $X$ be an $h$-homogeneous compact space, and let
$G=\Aut(X)$. Let $\Th=E_0(X)$ be the semigroup of all closed relations
$R$ on $X$ such that $\Dom R=\Ran R=X$.
We identify $G$ with the set of all invertible elements of $\Th$.
We prove in this section that $\Th$ can be identified
with the Roelcke compactification of $G$.

The space $\Th$, being compact,
has a unique compatible uniformity. Let $\sU$ be the uniformity that
$G$ has as a subspace of $\Th$. The first part of Theorem~\name{main}
is equivalent to the following:

\proclaim{\No{Roel}. Theorem}
Let $X$ be an $h$-homogeneous compact space, $\Th=E_0(X)$, 
and $G=\Aut(X)$. Identify $G$ with the set of all invertible
elements of $\Th$. Then:
\roster
\item $G$ is dense in $\Th$;
\item the uniformity $\sU$  induced by the embedding of $G$ into $\Th$
coincides with the Roelcke uniformity $\sL\wedge\sR$ on $G$.
\endroster
\endproclaim

Let us first introduce some notation.
Let $\g=\{U_\a:\a\in A\}$ be a finite clopen partition of $X$.
A {\it $\g$-rectangle\/} is a set of the form $U_\a\ti U_\beta$,
$\a,\beta\in A$.
Given a relation $R\in \Th$, denote by $M(\g,R)$ the set
of all pairs $(\a,\beta)\in A\ti A$ such that $R$ meets
the rectangle $U_\a\ti U_\beta$. 
Let $\sV(\g,R)$ be the family $\{U_\a\ti U_\beta:(\a,\beta)\in M(\g,R)\}$
of all $\g$-rectangles which meet $R$.
If $r$ is a subset 
of $A\ti A$, put
$$
O_{\g,r}=\{R\in \Th: M(\g,R)=r\}.
$$ 
The sets of the form $O_{\g,r}$ constitute a base of $\Th$.
Denote by $E_0(A)$ the set of all relations
$r$ on $A$ such that $\Dom r=\Ran r=A$.
A set $O_{\g,r}$ is non-empty if and only if $r\in E_0(A)$.

Let $O_\g(R)$ be the set of all relations $S\in \Th$ which meet the
same $\g$-rectangles as $R$. 
We have $O_\g(R)=O_{\g,r}$,
where $r=M(\g,R)$.
The sets of the form $O_\g(R)$ constitute a base at $R$. If
$\l$ is another clopen partition of $X$ which refines $\g$,
then $O_\l(R)\sbs O_\g(R)$.

\demo{Proof of Theorem~\name{Roel}} 
Our proof proceeds in three parts.

(a): We prove that $G$ is dense in $\Th$. 

Let $\g=\{U_\a:\a\in A\}$ be a finite clopen partition of $X$ and
$r\in E_0(A)$. 
We must prove that $O_{\g,r}$ meets $G$. Decomposing each $U_\a$ into
a suitable number of clopen pieces, we can find a clopen partition
$\{W_{\a,\beta}:(\a,\beta)\in r\}$ of $X$ such that 
$U_\a=\bigcup\{W_{\a,\beta}:(\a,\beta)\in r\}$ for every $\a\in A$.
Similarly, there exists a clopen partition 
$\{W_{\a,\beta}':(\a,\beta)\in r\}$ of $X$ such that 
$U_\beta=\bigcup\{W_{\a,\beta}':(\a,\beta)\in r\}$ for every $\beta\in A$.
Let $f\in G$ be a self-homeomorphism of $X$ such that $f(W_{\a,\beta})=
W_{\a,\beta}'$ for every $(\a,\beta)\in r$. The graph of $f$ meets
each rectangle of the form $W_{\a,\beta}\ti W_{\a,\beta}'$, 
$(\a,\beta)\in r$, and is contained in the union of such rectangles. 
It follows that $M(\g,f)=r$ and $f\in G\cap O_{\g,r}\ne\emptyset$.

(b): We prove that the uniformity $\sU$ is coarser than $\sL\wedge\sR$. 

This is a special case of the following:

\proclaim{\No{compunif}. Lemma}
For every compact space $K$ the map
$h\mapsto \G(h)$ from $\Aut(K)$ to $E_0(K)$
is $\sL\wedge\sR$-uniformly continuous. 
\endproclaim

\demo{Proof}
It suffices to prove that the map under consideration is $\sL$-uniformly
continuous and $\sR$-uniformly continuous. Let $d$ be a continuous
pseudometric on $K$. Let $d_2$ be the pseudometric on $K^2$ defined
by $d_2((x,y),(x',y'))=d(x,x')+d(y,y')$, and let $d_H$ be the corresponding
Hausdorff pseudometric on $E_0(K)$. If $R,S\in E_0(K)$ and $a>0$,
then $d_H(R,S)\le a$ if and only if each of the relations $R$ and $S$
is contained in the closed $a$-neighbourhood of the other with respect
to $d_2$. The pseudometrics of the form $d_H$ generate the uniformity
of $E_0(K)$. 

Let $d_s$ be the right-invariant pseudometric on $\Aut(K)$ defined
by $d_s(f,g)=$\newline
$\sup\{d(f(x),g(x)):x\in K\}$. The pseudometrics of the
form $d_s$ generate the right uniformity $\sR$ on $\Aut(K)$.
Since $d_H(\G(f),\G(g))\le d_s(f,g)$, the map $\G:\Aut(K)\to E_0(K)$
is $\sR$-uniformly continuous. For the left uniformity $\sL$ we can either
use a similar argument, or note that the involution on $\Aut(K)$
is an isomorphism between $\sL$ and $\sR$, and use the formula
$\G(f)=\G(f\obr)^*$ to reduce the case of $\sL$ to the case of $\sR$.
\qed\enddemo

(c): We prove that $\sU$ is finer than $\sL\wedge\sR$. 

Let $\g=\{U_\a:\a\in A\}$ be a finite clopen partition of $X$.
Put $V_\g=\{f\in G:f(U_\a)=U_\a\text{ for every }\a\in A\}$.
The open subgroups of the form $V_\g$ constitute a base at unity of $G$.
We must show that if $f,g\in G$ are close enough in $\Th$, then
$f\in V_\g g V_\g$.

The set of all pairs $(R,S)\in \Th^2$ such that $M(\g,R)=M(\g,S)$
is a neighbourhood of the diagonal in $\Th^2$ and therefore an entourage
for the unique compatible uniformity on $\Th$.
It suffices to prove that for every $f,g\in G$ the condition
$M(\g,f)=M(\g,g)$ implies that $f\in V_\g g V_\g$.
Suppose that $M(\g,f)=M(\g,g)=r$.
The following conditions are equivalent for every $\a,\beta\in A$:
(1) $f(U_\a)\cap U_\beta\ne\emptyset$; (2) $g(U_\a)\cap U_\beta\ne\emptyset$;
(3) $(\a,\beta)\in r$. Pick $u\in G$ such that 
$u(f(U_\a)\cap U_\beta)=g(U_\a)\cap U_\beta$ for every $(\a,\beta)\in r$.
Such a self-homeomorphism $u$ of $X$ exists, since all non-empty
clopen subsets
of $X$ are homeomorphic. Since for a fixed $\beta\in A$ the sets
$f(U_\a)\cap U_\beta$ cover $U_\beta$, we have $u(U_\beta)=U_\beta$.
Thus $u\in V_\g$. It follows that $uf(U_\a)\cap U_\beta=
u(f(U_\a)\cap U_\beta)=g(U_\a)\cap U_\beta$ for all $\a,\beta\in A$
and hence $uf(U_\a)=g(U_\a)$ for every $\a\in A$. Put $v=g\obr uf$.
Since $uf(U_\a)=g(U_\a)$, we have $v(U_\a)=U_\a$ for every $\a\in A$.
Thus $v\in V_\g$ and $f=u\obr gv\in V_\g g V_\g$.
\qed\enddemo

\newSection

\head
\S~\number\the\secno. Continuity-like properties of composition
\endhead

We preserve all the notation of the previous section. In particular,
$X$ is an $h$-homoge\-ne\-ous compact space, $G=\Aut(X)$, $\Th=E_0(X)$.

Recall that is a non-empty collection $\sF$ of non-empty subsets of a set $Y$
is a {\it filter base\/} on $Y$ if for every $A,B\in \sF$ there is $C\in \sF$
such that $C\sbs A\cap B$. If $Y$ is a topological space, $\sF$ is a filter
base on $Y$ and $x\in Y$, then $x$ is a {\it cluster point\/} of $\sF$ if
every neighbourhood of $x$ meets every member of $\sF$, and $\sF$ 
{\it converges\/} to $x$ if every neighbourhood of $x$ contains a member of 
$\sF$. If $\sF$ and $\sG$ are two filter bases on $G$, let
$\sF\sG=\{AB:A\in \sF,\, B\in \sG\}$. 

For every $R\in \Th$ let $\sF_R=\{G\cap V: V\text{ is a neighbourhood of $R$ in
$\Th$}\}$. In other words, $\sF_R$ is the trace on $G$ of the filter of
neighbourhoods of $R$ in $\Th$. We have noted that the multiplication
on $\Th$ is not continuous. 
If $R,S\in \Th$, it is not true in general that $\sF_R\sF_S$ converges to 
$RS$. However, $RS$ is a cluster point of $\sF_R\sF_S$. This fact
will be used in the next section.

\proclaim{\No{cluster}. Proposition}
If $R,S\in \Th$, then $RS$ is a cluster point of the filter base
$\sF_R\sF_S$.
\endproclaim

We need some lemmas. First we
note that for any compact space $K$ the composition
of relations is upper-semicontinuous on $E(K)$ in the following
sense:

\proclaim{\No{lemupper}. Lemma}
Let $K$ be a compact space, $R,S\in E(K)$. Let $O$ be an open
set in $K^2$ such that $RS\sbs O$. Then there exist open
sets $V_1, V_2$ in $K^2$ such that $R\sbs V_1$, $S\sbs V_2$,
and for every $R', S'\in E(K)$ such that $R'\sbs V_1$, $S'\sbs V_2$
we have $R'S'\sbs O$.
\endproclaim

\demo{Proof}
Consider the following three closed sets in $K^3$: 
$F_1=\{(x,z,y):(z,y)\in R\}$, $F_2=\{(x,z,y):(x,z)\in S\}$,
$F_3=\{(x,z,y):(x,y)\notin O\}$. The intersection of these three
sets is empty. There exist neighbourhoods of these sets with empty
intersection. We may assume that the neighbourhoods of $F_1$ and $F_2$
are of the form $\{(x,z,y):(z,y)\in V_1\}$ and $\{(x,z,y):(x,z)\in V_2\}$,
respectively, where $V_1$ and $V_2$ are open in $K^2$. The sets $V_1$ and 
$V_2$ are as required.
\qed\enddemo

\proclaim{\No{lemcomp}. Lemma}
Let $\g=\{U_\a:\a\in A\}$ be a finite clopen partition of $X$. For every
$R,S\in \Th$ we have $M(\g,RS)\sbs M(\g,R)M(\g,S)$ {\rm(the product on the right
means the composition of relations on $A$).}
\endproclaim

\demo{Proof}
Let $(\a,\beta)\in M(\g,RS)$. Then $RS$ meets the rectangle $U_\a\ti U_\beta$.
Pick $(x,y)\in RS\cap (U_\a\ti U_\beta)$. There exists $z\in X$ 
such that $(x,z)\in S$ and $(z,y)\in R$. Pick $\d\in A$ such that 
$z\in U_\d$. Then $(x,z)\in S\cap (U_\a\ti U_\d)$, 
$(z,y)\in R\cap (U_\d\ti U_\beta)$, hence $(\a,\d)\in M(\g,S)$ and 
$(\d,\beta)\in M(\g, R)$. It follows that $(\a,\beta)\in M(\g,R)M(\g,S)$.
\qed\enddemo

\proclaim{\No{lemappr}. Lemma}
Let $\l=\{U_\a:\a\in A\}$ be a finite clopen partition of $X$, and let
$r,s\in E_0(A)$. There exist $f,g\in G$ such that $M(\l,f)=r$,
$M(\l,g)$=s and $M(\l,fg)=rs$.
\endproclaim

\demo{Proof}
We modify the proof of Theorem~\name{Roel}. For every $\g\in A$
take a clopen partition $\{V_{\a,\g,\beta}:(\a,\g)\in s,\, (\g,\beta)\in r\}$
of $U_\g$. For every $(\g,\beta)\in r$ put $W_{\g,\beta}=\bigcup
\{V_{\a,\g,\beta}:(\a,\g)\in s\}$. For every $(\a,\g)\in s$ put 
$Y_{\a,\g}'=\bigcup
\{V_{\a,\g,\beta}:(\g,\beta)\in r\}$. Take a clopen partition
$\{W_{\g,\beta}':(\g,\beta)\in r\}$ of $X$ such that for every 
$\beta\in A$ we have $U_\beta=\bigcup\{W_{\g,\beta}':(\g,\beta)\in r\}$.
Take a clopen partition
$\{Y_{\a,\g}:(\a,\g)\in s\}$ of $X$ such that for every 
$\a\in A$ we have $U_\a=\bigcup\{Y_{\a,\g}:(\a,\g)\in s\}$.
There exist $f\in G$ such that $f(W_{\g,\beta})=W_{\g,\beta}'$ for
every $(\g,\beta)\in r$. There exists $g\in G$ such that 
$g(Y_{\a,\g})=Y_{\a,\g}'$ for every $(\a,\g)\in s$. The graph of $f$
meets every rectangle $W_{\g,\beta}\ti W_{\g,\beta}'$, $(\g,\beta)\in r$,
and is contained in the union of such rectangles. Since
$W_{\g,\beta}\ti W_{\g,\beta}'\sbs U_\g\ti U_\beta$, it follows that 
$M(\l,f)=r$. Similarly, $M(\l,g)=s$. We claim that $M(\l,fg)=rs$.
Let $(\a,\beta)\in rs$. There exists $\g\in A$ such that $(\a,\g)\in s$
and $(\g,\beta)\in r$. We have $g(U_\a)\supset g(Y_{\a,\g})=Y_{\a,\g}'\supset
V_{\a,\g,\beta}$ and $f\obr(U_\beta)\supset f\obr(W_{\g,\beta}')
=W_{\g,\beta}\supset V_{\a,\g,\beta}$. Thus $V_{\a,\g,\beta}\sbs g(U_\a)\cap
f\obr(U_\beta)\ne\emptyset$. It follows that the graph of $fg$ meets
the rectangle $U_\a\ti U_\beta$. This means that $(a,\beta)\in M(\l,fg)$.
We have proved that $rs\sbs M(\l,fg)$. The reverse inclusion follows
from Lemma~\name{lemcomp}.
\qed\enddemo

\demo{Proof of Proposition~\name{cluster}}
Let $U_1$, $U_2$, $U_3$ be neighbourhoods in $\Th$ of $R$, $S$ and $RS$,
respectively. We must show that $U_3$ meets the set 
$(U_1\cap G)(U_2\cap G)$.

Fix a clopen partition $\l$ of $X$ such that $O_\l(RS)\sbs U_3$.
Lemma~\name{lemupper} implies that there exists a clopen partition
$\g$ of $X$ such that for every $R'\in O_\g(R)$ and $S'\in O_\g(S)$
we have  $R'S'\sbs \bigcup \sV(\l,RS)$ (recall that $\sV(\l,RS)$
is the family of all $\l$-rectangles that meet $RS$).
We may assume that $\g$ refines $\l$ and that
$O_\g(R)\sbs U_1$, $O_\g(S)\sbs U_2$. Put $r=M(\g,R)$, $s=M(\g,S)$.
According to Lemma~\name{lemappr}, there exist $f,g\in G$ such that
$M(\g,f)=r$, $M(\g,g)=s$ and $M(\g,fg)=rs$. Then
$f\in G\cap O_\g(R)$ and $g\in G\cap O_\g(S)$.
Lemma~\name{lemcomp} implies that $M(\g,RS)\sbs rs=M(\g,fg)$.
This means that (the graph of)
$fg$ meets every member of the family $\sV(\g, RS)$. 
Then every member of $\sV(\l,RS)$ meets $fg$,
since every member of $\sV(\l,RS)$ contains a member of $\sV(\g,RS)$.
On the other hand, by the choice of $\g$ we have $fg\sbs \bigcup \sV(\l,RS)$.
It follows that $M(\l,fg)=M(\l,RS)$. Thus 
$fg\in O_\l(RS)\sbs U_3$ and hence $fg\in 
(U_1\cap G)(U_2\cap G)\cap U_3\ne\emptyset$.
\qed\enddemo

\newSection

\head
\S~\number\the\secno. Proof of Main Theorem, part 2
\endhead

Let $X$, as before, be a compact
$h$-homogeneous space, $G=\Aut(X)$, $\Th=E_0(X)$. 
We saw that $G$ is Roelcke-precompact and
that $\Th$ can be identified with the Roelcke compactification of $G$. 
In this section we prove that $G$ is minimal and topologically simple.

If $H$ is a group and $g\in H$, we denote by $l_g$ (respectively, $r_g$)
the left shift of $H$ defined by $l_g(h)=gh$ (respectively,
the right shift defined by $r_g(h)=hg$). 

\proclaim{\No{shiftroel}. Proposition}
Let $H$ be a topological group, and let $K$ be the Roelcke completion of $H$.
Let $g\in H$. 
Each of the following self-maps of $H$ extends to a self-homeomorphism of $K$:
(1)~the left shift~$l_g$; (2)~the right shift~$r_g$; (3)~the inversion 
$g\mapsto g\obr$. 
\endproclaim

\demo{Proof}
Let $\sL$ and $\sR$ be the left and the right uniformity on $H$, respectively.
In each of the cases (1)--(3) 
the map $f:H\to H$ under consideration is an automorphism of
the uniform space $(H,\sL\wedge\sR)$. This is obvious for the case~(3).
For the cases~(1) and~(2), observe that the uniformities $\sL$ and
$\sR$ are invariant under left and right shifts, hence the same is true for
their greatest lower bound $\sL\wedge\sR$. It follows that in all cases
$f$ extends to an automorphism of the completion $K$ of the uniform space
$(H,\sL\wedge\sR)$.
\qed\enddemo

For $g\in G$ define self-maps $L_g:\Th\to \Th$ and $R_g:\Th\to \Th$
by $L_g(R)=gR$ and $R_g(R)=Rg$. 

\proclaim{\No{contshift}. Proposition}
For every $g\in G$ the maps $L_g:\Th\to \Th$ and $R_g:\Th\to \Th$
are continuous.
\endproclaim

\demo{Proof}
We have $gR=\{(x,g(y)):(x,y)\in R\}$. 
Let $\l=\{U_\a:\a\in A\}$ be a clopen partition of $X$.
Let $r=M(\l,gR)$, and let
$O_\l(gR)=\{S\in \Th: M(\l,S)=r\}$ be a basic neighbourhood
of $gR$. Let $U$ be the set of all $T\in \Th$ such that 
$T$ meets every member of the family $\{U_\a\ti g\obr(U_\beta):
(\a,\beta)\in r\}$ and is contained in the union of this family.
Then $U$ is a neighbourhood of $R$ and $L_g(U)=O_\l(gR)$.
Thus $L_g$ is continuous. The argument for $R_g$ is similar.
\qed\enddemo

Let $\D$ be the diagonal in $X^2$. 

\proclaim{\No{maximidem}. Proposition}
Let $S$ be a closed subsemigroup of $\Th$, and let $T$ be the set of all 
$p\in S$ such that $p\supset\D$. If $T\ne\emptyset$,
then $T$ has a greatest element $p$, and $p$ is an idempotent.
\endproclaim


\demo{Proof}
We claim that every non-empty closed subset of $\Th$ has a maximal element.
Indeed, if $C$ is a non-empty
linearly ordered subset of $\Th$, 
then $C$ has a least upper bound $b=\cl{\cup C}$ in $\Th$,
and $b$ belongs to the closure of $C$ in $\Th$. Thus our claim follows from
Zorn's lemma. 

The set $T$ is a closed subsemigroup of $\Th$.
Let $p$ be a maximal element of $T$.
For every $q\in T$ we have $pq\supset p\D=p$, whence $pq=p$.
It follows that $p$ is an idempotent and that 
$p=pq\supset\D q=q$. Thus $p$ is the greatest element of $T$.
\qed\enddemo

An {\it inner automorphism\/} of $\Th$ is a map of the form 
$p\mapsto gpg\obr$, $g\in G$. 

\proclaim{\No{innidemp}. Proposition}
There are precisely two elements in $\Th$ 
which are invariant under all inner automorphisms of $\Th$, namely 
$\D$ and $X^2$.
\endproclaim

\demo{Proof}
A relation $R\in \Th$ is invariant under all inner automorphisms
if and only if the following holds: if $x,y\in X$ and $(x,y)\in R$,
then $(f(x),f(y))\in R$ for every $f\in G$. Suppose that $R\in \Th$ has
this property and $\D\ne R$. Pick $(x,y)\in R$ such that $x\ne y$.
We claim that the set $B=\{(f(x),f(y)):f\in G\}$ is dense in $X^2$. Indeed,
pick disjoint clopen neighbourhoods $U_1$ and $U_2$  of $x$ and $y$, 
respectively, such that $X$ is not covered by $U_1$ and $U_2$. Given
disjoint clopen non-empty sets $V_1$ and $V_2$, by $h$-homogeneity
of $X$ we can find
an $f\in G$ such that $f(U_i)\sbs V_i$, $i=1,2$.
It follows that $V_1\ti V_2$ meets $B$, hence $B$ is dense in $X^2$.
Since $B\sbs R$, it follows that $R=X^2$.
\qed\enddemo 

\proclaim{\No{compnorm}. Proposition}
The group $G$ has no compact normal subgroups other than $\{1\}$.
\endproclaim

We shall prove later that actually $G$ has no non-trivial closed normal 
subgroups.

\demo{Proof}
Let $H\ne\{1\}$ be a normal subgroup of $G$. We show that $H$ is not compact.

Let $Y$ be the collection of all non-empty
clopen sets in $X$. Consider $Y$ as
a discrete topological space. The group $G$ has a natural continuous
action on $Y$. Pick $f\in H$, $f\ne 1$. Pick $U\in Y$
such that $f(U)\cap U=\emptyset$ and $X\setminus(f(U)\cup U)\ne\emptyset$.
Let $Y_1$ be the set of all $V\in Y$ such that $V$
is a proper subset of $X\setminus U$. 
If $V\in Y_1$,
there exists $h\in G$ such that $h(U)=U$ and $h(f(U))=V$. Put $g=hfh\obr$.
Then $g(U)=V$.
Since $H$ is a normal subgroup of $G$, we have $g\in H$. It follows
that the $H$-orbit of $U$ contains $Y_1$. Since $Y_1$ is infinite,
$H$ cannot be compact.
\qed\enddemo

\proclaim{\No{simplemin}. Proposition}
For every topological group $H$ the following
conditions are equivalent:
\roster
\item $H$ is minimal and topologically simple;
\item if $f:H\to H'$ is a continuous onto homomorphism of topological groups,
then either $f$ is a homeomorphism, or $|H'|=1$.\qed
\endroster
\endproclaim

We are now ready to prove Theorem~\name{main}, part~(2):

\smallskip
{\it
For every compact $h$-homogeneous space $X$ 
the topological group $G=\Aut(X)$
is minimal and topologically simple.
}
\smallskip

\demo{Proof}
Let $f:G\to G'$ be a continuous onto homomorphism.
According to Proposition~\name{simplemin}, 
we must prove that either $f$ is a homeomorphism
or $|G'|=1$. 

Since $G$ is Roelcke-precompact, so is $G'$. Let $\Th'$ be the Roelcke
compactification of $G'$. The homomorphism $f$ extends to a continuous map
$F:\Th\to\Th'$. Let $e'$ be the unity of $G'$, and let $S=F\obr(e')\sbs \Th$.

Claim 1. $S$ is a subsemigroup of $\Th$.

Let $p,q\in S$. In virtue of Proposition~\name{cluster}, there exist filter
bases $\sF_p$ and $\sF_q$ on $G$ such that $\sF_p$ converges to $p$ (in $\Th$),
$\sF_q$ converges to $q$ and $pq$ is a cluster point of the filter base
$\sF_p\sF_q$. The filter bases $\sF_p'=F(\sF_p)$ and $\sF_q'=F(\sF_q)$ on $G'$
converge to $F(p)=F(q)=e'$, hence the same is true for the filter base
$\sF_p'\sF_q'=F(\sF_p\sF_q)$. Since $pq$ is a cluster point of
$\sF_p\sF_q$, $F(pq)$ is a cluster point of the convergent filter base
$F(\sF_p\sF_q)$. A convergent filter on a Hausdorff space has only one cluster
point, namely the limit. Thus $F(pq)=e'$ and hence $pq\in S$.

Claim 2. The semigroup $S$ is closed under involution.

In virtue of Proposition~\name{shiftroel},
the inversion on $G'$ extends to an involution $x\mapsto x^*$ of $\Th'$. Since
$F(p^*)=F(p)^*$ for every $p\in G$, the same holds for every $p\in \Th$. Let
$p\in S$. Then $F(p^*)=F(p)^*=e'$ and hence $p^*\in S$.

Claim 3. If $g\in G$ and $g'=f(g)$, then $F\obr(g')=gS=Sg$.

We saw that the left shift $h\mapsto gh$ of $G$ extends to a continuous
self-map $L=L_g$ of $\Th$ defined by $L(p)=gp$
(Proposition~\name{contshift}). 
According to Proposition~\name{shiftroel}, the self-map $x\mapsto g'x$ of $G'$
extends to self-homeomorphism $L'$ of $\Th'$. The maps $FL$ and $L'F$ from
$\Th$ to $\Th'$ coincide on $G$ and hence everywhere. Replacing
$g$ by $g\obr$, we see that $FL\obr=(L')\obr F$.
Thus
$F\obr(g')=F\obr L'(e')=LF\obr(e')=gS$. Using right shifts instead of
left shifts, we similarly conclude that $F\obr(g')=Sg$.

Claim 4. $S$ is invariant under inner automorphisms of $\Th$.

We have just seen that $gS=Sg$ for every $g\in G$, hence
$gSg\obr=S$.

Let $T=\{r\in S:r\supset \D\}$. According to Proposition~\name{maximidem},
there is a greatest element $p$ in $T$.
Claim~4
implies that $p$ is invariant under inner automorphisms.
In virtue of Proposition~\name{innidemp}, either $p=\D$ or $p=X^2$.
We shall show that either $f$ is a homeomorphism or $|G'|=1$,
according to which of the cases $p=\D$ or $p=X^2$ holds.

First assume that $p=\D$.

Claim 5 ($p=\D$). All elements of $S$ are invertible in $\Th$.

Let $r\in S$. Then $r^*r\in S$
and $rr^*\in S$, since $S$ is a symmetrical semigroup.
Since $\Dom r=\Ran r=X$, we have 
$r^*r\supset \D$ and $rr^*\supset\D$.
The assumption $p=\D$ implies that every relation $s\in S$ 
such that $s\supset\D$ must be equal to $\D$.
Thus $rr^*=r^*r=\D$ and $r$ is invertible.

Claim 6 ($p=\D$). $|S|=1$.

Claim 5 implies that $S$ is a subgroup of $G$.
This subgroup is normal (Claim~4) and compact, since $S$ is closed in $\Th$.
Proposition~\name{compnorm} implies that $|S|=1$.

Claim 7 ($p=\D$). $f:G\to G'$ is a homeomorphism.

Claims~6 and~3 imply that $G= F\obr(G')$ and that the map $f:G\to G'$ is
bijective. Since $F$ is a map between compact spaces, it is perfect, and hence
so is the map $f:G=F\obr(G')\to G'$. Thus $f$, being a perfect bijection,
is a homeomorphism.

Now consider the case $p=X^2$.

Claim 8. If $p=X^2\in S$, then $G'=\{e'\}$.

Let $g\in G$ and $g'=f(g)$. We have $gp=p\in S$. On the other hand,
Claim~3 implies that $gp\in gS=F\obr(g')$. Thus $g'=F(gp)=
F(p)=e'$.
\qed\enddemo

\newSection

\head
\S~\number\the\secno. Remarks
\endhead

The group $\Aut(K)$ is Roelcke-precompact also for some compact spaces $K$
which are not zero-dimensional. For example, let $I=[0,1]$ and $G=\Aut(I)$. 
Identify
$G$ with a subspace of $E(I)$, as above. The Roelcke compactification
of $G$ can be identified with the closure of $G$ in $E(I)$. Let $G_0$
be the subgroup of all $f\in G$ which leave the end-points of the interval 
$I$ fixed. The closure of $G_0$ in $E(I)$ is the set of all curves $c$ in
the square $I^2$ such that $c$ connects the points $(0,0)$ and $(1,1)$
and has the following property: there are no points $(x,y)\in c$ and
$(x',y')\in c$ such that $x<x'$ and $y>y'$. This can be used to yield
an alternative proof of D.~Gamarnik's theorem saying that $G$ is minimal
[\bibl{gamarn}].

Let $K=I^\o$ be the
Hilbert cube and $G=\Aut(K)$. I do not know if $G$ is minimal or
Roelcke-precompact in this case.

\enddocument